\documentclass{article}
\usepackage[letterpaper,top=2cm,bottom=2cm,left=2.5cm,right=2.5cm,marginparwidth=1.75cm]{geometry}

\usepackage{mathtools}
\usepackage{amsmath}
\usepackage{amssymb}
\usepackage{amsthm}
\usepackage{latexsym}
\usepackage[square,numbers]{natbib}
\usepackage{authblk}
\usepackage{soul}

\numberwithin{equation}{section}

\author{Efe Gürel$^{1}$\\
        \small $^{1}$TÜBİTAK Natural Sciences High School, Kocaeli 41400, Turkey    }
\date{}

\newtheorem{theorem}{Theorem}[section]

\newtheorem{corollary}[theorem]{Corollary}

\newtheorem{lemma}[theorem]{Lemma}

\theoremstyle{definition}
\newtheorem{example}{Example}

\title{A Recipe For Obtaining Algebraic Addition Theorems Of The Weierstrass Elliptic Function}

\begin{document}

\maketitle
\begin{abstract}
    In this paper, we present a general method for obtaining addition theorems of the Weierstrass elliptic function $\wp(z)$ in terms of given parameters. We obtain the classical addition theorem for the Weierstrass elliptic function as a special case. Furthermore, we give novel two-term addition, three-term addition, duplication and triplication formulas. New identities for elliptic invariants are also proven.\\



\end{abstract}

\section{Introduction}
Let $\Omega=(2\omega_1,2\omega_2)\in \mathbb{C}^2$ be such that $\text{Im}(\omega_2/\omega_1)>0$ and $\Lambda_\Omega=2\omega_1\mathbb{Z}+2\omega_2\mathbb{Z}$. The Weierstrass elliptic function $\wp$ is defined as
\begin{align*}
    \wp(z)=\wp(z;\Lambda_\Omega)=\frac{1}{z^2}+\sum_{\omega\in\Lambda_\Omega^*}\frac{1}{(z-\omega)^2}-\frac{1}{\omega^2},
\end{align*}
where $\Lambda_\Omega^*=\Lambda_\Omega-\{0\}$ and the series converges absolutely and uniformly except at the poles $z\in\Lambda_\Omega$. It is an even elliptic function of order $2$ with double poles at $z\in\Lambda_\Omega$. Weierstrassian elliptic invariants of the lattice $\Lambda_\Omega$ are defined as
\begin{align*}
    g_2=60\sum_{\omega\in\Lambda_\Omega^*}\omega^{-4} \qquad \text{and}\qquad g_3=140\sum_{\omega\in\Lambda_\Omega^*}\omega^{-6}.
\end{align*}
The function $\wp$ satisfies the well-known algebraic differential equation
\begin{align}\label{WDE}
    \left( \wp'(z) \right)^2=4\wp(z)^3-g_2\wp(z)-g_3=4(\wp(z)-e_1)(\wp(z)-e_2)(\wp(z)-e_3).
\end{align}
Here $e_1=\wp(\omega_1)$, $e_2=\wp(\omega_2)$ and $e_3=\wp(\omega_1+\omega_2)$ are the half-period values. We have the famous formulae
\begin{gather*}
    e_1+e_2+e_3=0\\
    e_1e_2+e_2e_3+e_3e_1=-g_2/4\\
    e_1e_2e_3=g_3/4.
\end{gather*}
Let $\sigma(z;\Lambda_\Omega)$ denote the Weierstrass $\sigma$-function, defined by
\begin{align*}                                          \sigma(z)=\sigma(z;\Lambda_{\Omega})=z\prod_{\omega\in\Lambda_\Omega^*}\left( 1-\frac{z}{\omega} \right)\exp\left( \frac{z}{\omega}+\frac{z^2}{2\omega^2} \right).
\end{align*}
We have the basic relation $\left(\log\sigma(z)\right)''=-\wp(z)$. Algebraic addition theorems are fundamental objects in the theory of elliptic functions. Frobenius-Stickelberger determinant identity or pseudo-addition theorem \cite{Frobenius,WhittakerWatson} states that
\begin{align}\label{FrobeniusStick}
    \begin{vmatrix}
        1 & \wp(z_0) & \wp'(z_0) &\cdots &\wp^{(n-1)}(z_0) \\
        1 & \wp(z_1) & \wp'(z_1) &\cdots &\wp^{(n-1)}(z_1)\\
        \vdots & \vdots & \vdots &\vdots &\vdots \\
        1 & \wp(z_n) & \wp'(z_n) &\cdots & \wp^{(n-1)}(z_n)
    \end{vmatrix}=(-1)^{\frac{n(n-1)}{2}}1!2!\cdots n!\frac{\sigma(z_0+z_1+\cdots+z_n)}{\sigma(z_0)^{n+1}\cdots\sigma(z_n)^{n+1}}\prod_{0\le i < j \le n } \sigma(z_i-z_j),
\end{align}
and for $n=1$ reduces to the well-known formula  
\begin{align*}
    \wp(z)-\wp(w)=\frac{\sigma(w+z)\sigma(w-z)}{\sigma(w)^2\sigma(z)^2}.
\end{align*}
Frobenius-Stickelberger identities may be used to obtain algebraic addition theorems for the function $\wp$ (for example, see \cite{Chandra,Lawdens}). For more related determinant identities we refer the reader to \cite{Frobenius,WhittakerWatson,Hermite,Kiepert,Kiepert2}. The classical addition theorem for $\wp$ is given as follows.
\begin{theorem}\label{ClassicAddition}
    Let $z,w,z\pm w\notin\Lambda_\Omega$, then we have
    \begin{align*}
        \begin{vmatrix}
            1 & \wp(z) & \wp'(z)  \\
            1 & \wp(w) & \wp'(w)  \\
            1 & \wp(z+w) & -\wp'(z+w) 
        \end{vmatrix}=0
    \end{align*}
    Furthermore,
    \begin{align*}
        \wp(z+w)=\frac{1}{4}\left( \frac{\wp'(z)-\wp'(w)}{\wp(z)-\wp(w)} \right)^2-\wp(z)-\wp(w).
    \end{align*}
\end{theorem}
Upon taking the limit $w\to z$, we also obtain the duplication formula
\begin{align*}
    \wp(2z)=\frac{1}{4}\left( \frac{\wp''(z)}{\wp'(z)} \right)^2-2\wp(z).
\end{align*}
The addition formula for the $\wp$ function is crucial and in fact a characteristic property of the $\wp$ function, as shown in \cite{DeterminantAddition}. For recent advances in addition theorems for elliptic and abelian functions, see \cite{EilbeckAbelian,EilbeckWeierstrass}. Recent work in the research of addition formulas by these authors uses a more algebraic geometric approach and is applicable to the general case of higher genus curves. The authors have a long series of publications discussing such investigations. However, the addition formulas obtained are not explicit in the classical sense. For example, the general addition theorems obtained for the classical elliptic (genus 1) case make use of conjugate variables. Even in the Weierstrass case, calculation of conjugate variables reduce to solving the equation $\wp'(u)=\wp'(v)=\wp'(w)$, which can not be easily done. The conjugate variables can be found explicitly only in the cases of special lattices which have particular symmetries, such as the equianharmonic and the lemniscatic case. Generally, one obtains Frobenius-Stickelberger like determinant identities. In this note, we discuss a recipe to obtain algebraic addition theorems for $\wp$, expressing $\wp(z_1+\cdots+z_\ell)$ in terms of $\wp^{(j)}(z_i)$ for predetermined $j$. The novelty of our approach is the immense flexibility in ones ability to choose parameters. One can choose which derivatives of the function $\wp$ and what constants can appear and modify them to obtain different addition formulas. Furthermore, our equations make no use of conjugate variables and only depend on the values of given variables $z_1,\ldots,z_\ell$. Another thing to note is that we make no assumptions on the lattice.

\section{Main Results}
For convenience, let $\wp^{(-2)}(z)\equiv1$ denote the constant function. Thus $\wp^{(i)}$ is an elliptic function of order $i+2$. Throughout this chapter, we assume that $\{n_i\}_{i=1}^m$,  $\{k_i\}_{i=1}^\ell$ are integers such that $n_i\in \mathbb{N}\cup \{-2\}$ for every $i=1,\cdots,m$ and $k_i\in \mathbb{N}\cup \{-2\}$ for every $i=1,\cdots,\ell$. Assume $n_i\neq n_j$ for $i\neq j$ and $k_i\neq k_j$ for $i\neq j$. Let $\{\gamma_i\}_{i=1}^{m}$ and $\{z_i\}_{i=1}^{\ell}$ be given complex numbers such that $z_i\notin\Lambda_\Omega$ for every $i=1,\ldots,\ell$ and
\begin{align}\label{DetCondition}
    \begin{vmatrix}
        \wp^{(k_1)}(z_1) & \wp^{(k_2)}(z_2) &\cdots &\wp^{(k_\ell)}(z_1) \\
        \wp^{(k_1)}(z_2) & \wp^{(k_2)}(z_2) &\cdots &\wp^{(k_\ell)}(z_2)\\
        \vdots & \vdots & \vdots &\vdots  \\
        \wp^{(k_1)}(z_\ell) & \wp^{(k_2)}(z_\ell) &\cdots & \wp^{(k_\ell)}(z_\ell)
    \end{vmatrix}\neq0.
\end{align}
Note that this implies $z_i\not\equiv z_j \bmod \Lambda_\Omega$ for $i\neq j$. Let $n_a=\max n_i$, $k_b=\max k_i$ and $\ell=\max(n_a,k_b)+1$. Thus there exists complex numbers $\{\lambda_i\}_{i=1}^{\ell}$ such that the system
\begin{align}\label{ZJSystem}
    \sum_{i=1}^{m}\gamma_i \wp^{(n_i)}(z_j)=\sum_{i=1}^{\ell}\lambda_i\wp^{(k_i)}(z_j)
\end{align}
holds for every $j=1,\ldots, \ell$. Let us denote $z=z_1+\cdots+ z_\ell$. Finally, assume that $z\notin \Lambda_\Omega$ and both of the equations $n_a=k_b$, $\gamma_a=\lambda_b$ do not hold at the same time. In order to obtain addition formulas, we shall prove that $-z$ is another solution of the equation \eqref{ZJSystem}. The assumptions above guarantee that the poles of the function $\wp$ are not encountered. We will see in the future that these assumptions may be relaxed.
\newline

We first prove the determinant form of the addition theorem.

\begin{theorem}\label{DetThm}
    The following equation holds,
    \begin{align*}
        \begin{vmatrix}
        \sum_{i=1}^{m}\gamma_i \wp^{(n_i)}(z_1)& \wp^{(k_1)}(z_1) &\cdots &\wp^{(k_\ell)}(z_1) \\
        \sum_{i=1}^{m}\gamma_i \wp^{(n_i)}(z_2)& \wp^{(k_1)}(z_2)  &\cdots &\wp^{(k_\ell)}(z_2)\\
        \vdots& \vdots  & \vdots &\vdots  \\
        \sum_{i=1}^{m}\gamma_i \wp^{(n_i)}(z_\ell)& \wp^{(k_1)}(z_\ell)  &\cdots & \wp^{(k_\ell)}(z_\ell) \\
        \sum_{i=1}^{m}\gamma_i \wp^{(n_i)}(-z) & \wp^{(k_1)}(-z)  &\cdots & \wp^{(k_\ell)}(-z) 
    \end{vmatrix}=0.
    \end{align*}
\end{theorem}
\begin{proof}
    Let us consider the function
    \begin{align*}
        \psi_{n,k,\gamma,\lambda}(s)=\sum_{i=1}^{m}\gamma_i \wp^{(n_i)}(s)-\sum_{i=1}^{\ell}\lambda_i\wp^{(k_i)}(s).
    \end{align*}
    By our assumptions, $\psi_{n,k,\gamma,\lambda}$ is an elliptic function of order $\max(n_a,k_b)+2=\ell+1$ with a pole of order $\ell+1$ at $s=0$. Since the equation \eqref{ZJSystem} holds for any $j=1,\ldots,\ell$, the function $\psi_{n,k,\gamma,\lambda}$ has $\ell$ distinct zeros, namely $z_1,\ldots, z_\ell$. Let $w$ denote the $\ell+1$th zero. Using Abel's theorem we obtain $w+\sum_{i=1}^{\ell}z_i\equiv0\bmod\Lambda_\Omega$ and $ w\equiv-z \bmod\Lambda_\Omega.$ Thus $-z$ is a root of $\psi_{n,k,\gamma,\lambda}$ and we get
    \begin{align}\label{PsiHasRoot-z}
        \sum_{i=1}^{m}\gamma_i \wp^{(n_i)}(-z)=\sum_{i=1}^{\ell}\lambda_i\wp^{(k_i)}(-z).
    \end{align}
    Putting equations \eqref{ZJSystem} and \eqref{PsiHasRoot-z} together gives
    \begin{align}\label{BigMatrixEq}
        \begin{pmatrix}
        \sum_{i=1}^{m}\gamma_i \wp^{(n_i)}(z_1)& \wp^{(k_1)}(z_1) &\cdots &\wp^{(k_\ell)}(z_1) \\
        \sum_{i=1}^{m}\gamma_i \wp^{(n_i)}(z_2)& \wp^{(k_1)}(z_2)  &\cdots &\wp^{(k_\ell)}(z_2)\\
        \vdots& \vdots  & \vdots &\vdots  \\
        \sum_{i=1}^{m}\gamma_i \wp^{(n_i)}(z_\ell)& \wp^{(k_1)}(z_\ell)  &\cdots & \wp^{(k_\ell)}(z_\ell) \\
        \sum_{i=1}^{m}\gamma_i \wp^{(n_i)}(-z) & \wp^{(k_1)}(-z)  &\cdots & \wp^{(k_\ell)}(-z) 
        \end{pmatrix}
        \begin{pmatrix}
            -1 \\
            \lambda_1\\
            \vdots\\
            \lambda_\ell
        \end{pmatrix}
        =-\begin{pmatrix}
            \psi_{n,k,\gamma,\lambda}(z_1) \\
            \psi_{n,k,\gamma,\lambda}(z_2)\\
            \vdots\\
            \psi_{n,k,\gamma,\lambda}(z_\ell)\\
            \psi_{n,k,\gamma,\lambda}(-z)
        \end{pmatrix}=0.
    \end{align}
    Since the vector on the left-hand side of \eqref{BigMatrixEq} is non-zero, the matrix on the left-hand side is singular. Therefore we have
    \begin{align*}
        \begin{vmatrix}
        \sum_{i=1}^{m}\gamma_i \wp^{(n_i)}(z_1)& \wp^{(k_1)}(z_1) &\cdots &\wp^{(k_\ell)}(z_1) \\
        \sum_{i=1}^{m}\gamma_i \wp^{(n_i)}(z_2)& \wp^{(k_1)}(z_2)  &\cdots &\wp^{(k_\ell)}(z_2)\\
        \vdots& \vdots  & \vdots &\vdots  \\
        \sum_{i=1}^{m}\gamma_i \wp^{(n_i)}(z_\ell)& \wp^{(k_1)}(z_\ell)  &\cdots & \wp^{(k_\ell)}(z_\ell) \\
        \sum_{i=1}^{m}\gamma_i \wp^{(n_i)}(-z) & \wp^{(k_1)}(-z)  &\cdots & \wp^{(k_\ell)}(-z) 
    \end{vmatrix}=0.
    \end{align*}
    This completes the proof.
\end{proof}
We remark that since $\wp^{(i)}(-z)=(-1)^i\wp^{(i)}(z)$ and $\wp$,$\wp^{(i)}$ satisfy an algebraic relationship for every $i$, Theorem \ref{DetThm} is indeed an algebraic addition theorem.
\newline

We need the following two lemmas.
\begin{lemma}\label{PDiffLemma}
    For all $n\in \mathbb{N}$, $\wp^{(2n)}\in\mathbb{Q}[g_2,g_2][\wp]$ and $\wp^{(2n+1)}\in\wp' \mathbb{Q}[g_2,g_2][\wp]$. Furthermore $\wp^{(2n)}$ is a polynomial in $\wp$ of degree $n+1$ and $\wp^{(2n+1)}/\wp'$  a polynomial in $\wp$ of degree $n$.
\end{lemma}
\begin{proof}
    The proof follows immediately after induction and equation \eqref{WDE}.
\end{proof}
\begin{lemma}\label{PhiLemma}
    The function defined as
    \begin{align}\label{PhiDef}
        \varphi_{n,k,\gamma,\lambda}(s)=\left( \sum_{2\nmid n_i}\gamma_i\wp^{(n_i)}(s)- \sum_{2\nmid k_i}\lambda_i\wp^{(k_i)}(s)\right)^2-\left( \sum_{2\mid n_i}\gamma_i\wp^{(n_i)}(s)- \sum_{2\mid k_i}\lambda_i\wp^{(k_i)}(s)\right)^2
    \end{align}
    is a polynomial in $\wp(s)$ of degree at most $\ell+1$.
\end{lemma}
\begin{proof}
    By Lemma \ref{PDiffLemma} and equation \eqref{WDE}, the first sum in \eqref{PhiDef} is in $\left(\wp' \mathbb{C}[\wp]\right)^2$ and also in $\mathbb{C}[\wp]$. Similarly the second sum in \eqref{PhiDef} is a polynomial in $\mathbb{C}[\wp]$. Therefore $\varphi_{n,k,\gamma,\lambda}\in\mathbb{C}[\wp]$. Again by Lemma \ref{PDiffLemma} and \eqref{WDE}, the first sum is a polynomial in $\wp$ of degree at most $\max_{2\nmid n_i,k_j}(n_i,k_j)+2$ and second sum is a polynomial in $\wp$ of degree at most $\max_{2\mid n_i,k_j}(n_i,k_j)+2$. Thus $\varphi_{n,k,\gamma,\lambda}$ is a polynomial in $\wp$ of degree at most $\max_{i,j}(n_i,k_j)+2=\ell+1$.
\end{proof}
We now see that there exists constants $\mu_{n,k,\gamma,\lambda}(r)\in\mathbb{C}$, $r=0,1,\ldots,\ell+1$ such that
\begin{align*}
    \varphi_{n,k,\gamma,\lambda}(s)=\sum_{r=0}^{\ell+1}\mu_{n,k,\gamma,\lambda}(r)\wp(s)^r.
\end{align*}
Let $S_r(x_1,\cdots x_n)$, $r\le n$ denote the elementary symmetric polynomials defined by
\begin{align*}
    S_{r}(x_1,\ldots,x_n)=\sum_{1\le i_1<\cdots< i_r\le n} x_{i_1} \ldots x_{i_r}.
\end{align*}
Furthermore, apply the convenience of notation
\begin{align*}
    S_{r}\wp(x_1,\ldots,x_n)=S_r\left( \wp(x_1),\ldots,\wp(x_n) \right).
\end{align*}
Assuming that $\varphi_{n,k,\gamma,\lambda}\not\equiv 0$, we now prove the explicit form of the addition theorem.
\begin{theorem}\label{MainAdditionThm}
    Let $z\not\equiv\pm z_i \bmod \Lambda_\Omega$ and $z_i\not\equiv- z_j \bmod \Lambda_\Omega$ for $i\neq j$. Then $\varphi_{n,k,\gamma,\lambda}(s)$ is a polynomial in $\wp(s)$ of degree $\ell+1$ with the distinct roots $\wp(z_1),\ldots,\wp(z_\ell),\wp(z)$. Furthermore,
    \begin{align*}
        S_{r}\wp(z_1,\ldots,z_\ell,z)=(-1)^r \frac{\mu_{n,k,\gamma,\lambda}(\ell+1-r)}{\mu_{n,k,\gamma,\lambda}(\ell+1)}.
    \end{align*}
\end{theorem}
\begin{proof}
    Rearranging equation \eqref{ZJSystem}, we get
    \begin{align*}
        \sum_{2\nmid n_i}\gamma_i\wp^{(n_i)}(z_j)- \sum_{2\nmid k_i}\lambda_i\wp^{(k_i)}(z_j)=- \sum_{2\mid n_i}\gamma_i\wp^{(n_i)}(z_j)+ \sum_{2\mid k_i}\lambda_i\wp^{(k_i)}(z_j).
    \end{align*}
    Therefore $\wp(z_1),\ldots,\wp(z_\ell)$ are distinct roots of $\varphi_{n,k,\gamma,\lambda}$. Using the fact that $\wp^{(i)}(-z)=(-1)^i\wp^{(i)}(z)$ and rearranging \eqref{PsiHasRoot-z}, we get
    \begin{align*}
        \sum_{2\nmid n_i}\gamma_i\wp^{(n_i)}(z)- \sum_{2\nmid k_i}\lambda_i\wp^{(k_i)}(z)= \sum_{2\mid n_i}\gamma_i\wp^{(n_i)}(z)- \sum_{2\mid k_i}\lambda_i\wp^{(k_i)}(z).
    \end{align*}
    Thus $\wp(z)$ is a root of $\varphi_{n,k,\gamma,\lambda}$. By Lemma \ref{PhiLemma} and assumption that $\varphi_{n,k,\gamma,\lambda}\not\equiv 0$, this proves that $\varphi_{n,k,\gamma,\lambda}(s)$ is a polynomial in $\wp(s)$ of degree $\ell+1$ with the distinct roots $\wp(z_1),\ldots,\wp(z_\ell),\wp(z)$. In particular, we see  $\mu_{n,k,\gamma,\lambda}(\ell+1)\neq0$. Therefore we have the factorization
    \begin{align*}
        \varphi_{n,k,\gamma,\lambda}(s)&=\mu_{n,k,\gamma,\lambda}(\ell+1)\left( \wp(s)-\wp(z) \right)\prod_{i=1}^{\ell}\left(\wp(s)-\wp(z_i)\right)\\
        &=\sum_{r=0}^{\ell+1}\mu_{n,k,\gamma,\lambda}(r)\wp(s)^r.
    \end{align*}    
    Comparing the coefficients of $\wp(s)^r$, we obtain
    \begin{align*}
        (-1)^r \mu_{n,k,\gamma,\lambda}(\ell+1)S_{r}\wp(z_1,\ldots,z_\ell,z)= \mu_{n,k,\gamma,\lambda}(\ell+1-r).
    \end{align*}
    Thus the proof is complete.
\end{proof}
\begin{theorem}\label{P(z)Explicit}
    Assume that $S_{r-1}\wp(z_1,\ldots,z_\ell)\neq 0$. Then we have
    \begin{align*}
        \wp(z)=\frac{1}{S_{r-1}\wp(z_1,\ldots,z_\ell)}\left((-1)^r \frac{\mu_{n,k,\gamma,\lambda}(\ell+1-r)}{\mu_{n,k,\gamma,\lambda}(\ell+1)} -S_{r}\wp(z_1,\ldots,z_\ell) \right).
    \end{align*}
\end{theorem}
\begin{proof}
    The proof follows from trivially by Theorem \ref{MainAdditionThm} and the identity
    \begin{align*}
        S_{r}\wp(z_1,\ldots,z_\ell,z)=\wp(z)S_{r-1}\wp(z_1,\ldots,z_\ell)+S_{r}\wp(z_1,\ldots,z_\ell).
    \end{align*}
\end{proof}
Taking $z_1,\ldots, z_\ell\to s$ in Theorem \ref{P(z)Explicit}, we obtain a multiplication theorem for $\wp(s)$. Comparing Theorem \ref{P(z)Explicit} for different $r_1,r_2$, we get the following corollary.
\begin{corollary}
      Let us denote $S_r=S_r\wp(z_1,\ldots,z_\ell)$. Then we have
     \begin{align*}
        (-1)^{r_1}S_{r_2-1}\mu_{n,k,\gamma,\lambda}(\ell+1-r_1)-(-1)^{r_2}S_{r_1-1}\mu_{n,k,\gamma,\lambda}(\ell+1-r_2)=\left( S_{r_1}S_{r_2-1}-S_{r_2} S_{r_1-1} \right)\mu_{n,k,\gamma,\lambda}(\ell+1).
     \end{align*}
\end{corollary}
We note that $\mu_{n,k,\gamma,\lambda}(r)$ are computable in terms of $\gamma_i, \lambda_j$ by Lemma \ref{PDiffLemma} since the expansion of $\wp^{(i)}$ in terms of $\wp$ is easily obtained by induction and differentiation. Applying  Cramer's rule to the system \eqref{ZJSystem} yields
\begin{align*}
    \lambda_j=\frac{ \begin{vmatrix}
        \wp^{(k_1)}(z_1) & \cdots&\sum_{i=1}^{m}\gamma_i \wp^{(n_i)}(z_1) &\cdots &\wp^{(k_\ell)}(z_1) \\
        \wp^{(k_1)}(z_2) & \cdots &\sum_{i=1}^{m}\gamma_i \wp^{(n_i)}(z_2) &\cdots &\wp^{(k_\ell)}(z_2)\\
        \vdots & \vdots & \vdots &\vdots &\vdots  \\
        \wp^{(k_1)}(z_\ell) & \cdots& \sum_{i=1}^{m}\gamma_i \wp^{(n_i)}(z_\ell) &\cdots & \wp^{(k_\ell)}(z_\ell)
    \end{vmatrix}}{ \begin{vmatrix}
        \wp^{(k_1)}(z_1) & \wp^{(k_2)}(z_2) &\cdots &\wp^{(k_\ell)}(z_1) \\
        \wp^{(k_1)}(z_2) & \wp^{(k_2)}(z_2) &\cdots &\wp^{(k_\ell)}(z_2)\\
        \vdots & \vdots & \vdots &\vdots  \\
        \wp^{(k_1)}(z_\ell) & \wp^{(k_2)}(z_\ell) &\cdots & \wp^{(k_\ell)}(z_\ell)
    \end{vmatrix}},
\end{align*}
where only $j$th column of the upper matrix has been replaced. We now remark that, by meromorphic continuation of both sides in above equations, one might drop the initial assumptions. The formulas given hold for all values of the variables appearing in them.

\section{Applications and Discussion}

In this section we illustrate our method with two examples, reobtaining the classical addition formula and novel three-term addition formulas. The applications are not limited to these examples. Here we restrict ourselves to the case of two and three variables because of the difficulty in calculations.

\begin{example}
    Let us take $n_1=\gamma_1=1$, $(k_1,k_2)=(0,-2)$ and $\ell=2$. We reobtain the classical proof of addition formula Theorem \ref{ClassicAddition}, as given in \cite{WhittakerWatson} and obtain two novel addition and duplication formulas in addition to new identities for elliptic invariants. We obtain the system
\begin{gather*}
    \wp'(z)=\lambda_1\wp(z)+\lambda_2\\
    \wp'(w)=\lambda_1\wp(w)+\lambda_2,
\end{gather*}
subject to the condition $\wp(z)\neq \wp(w)$. We then have $\psi(s)=\wp'(s)-\lambda_1\wp(s)-\lambda_2$ and $\varphi(s)=\left(\wp'(s)\right)^2-(\lambda_1\wp(s)+\lambda_2)^2$. By equation \eqref{WDE}, we calculate $\mu(r)$ as follows,
\begin{align*}
    &\mu(0)=-\lambda_2^2-g_3=-\left( \frac{\wp(z)\wp'(w)-\wp(w)\wp'(z)}{\wp(z)-\wp(w)} \right)^2-g_3\\
    &\mu(1)=-2\lambda_1\lambda_2-g_2=-2\frac{\left( \wp'(z)-\wp'(w) \right)\left( \wp(z)\wp'(w)-\wp(w)\wp'(z) \right) }{\left( \wp(z)-\wp(w) \right)^2}-g_2\\
    &\mu(2)=-\lambda_1^2=-\left( \frac{\wp'(z)-\wp'(w)}{\wp(z)-\wp(w)} \right)^2\\
    &\mu(3)=4.
\end{align*}
Thus we get the following formulae.
\begin{corollary}\label{ClassicalCorollary}
    We have the following addition formulae,
    \begin{align}\label{ClassicAdditionEq}
        \wp(z+w)&=\frac{1}{4}\left( \frac{\wp'(z)-\wp'(w)}{\wp(z)-\wp(w)} \right)^2-\wp(z)-\wp(w)\\
        \nonumber &=-\frac{\left( \wp'(z)-\wp'(w) \right)\left( \wp(z)\wp'(w)-\wp(w)\wp'(z) \right) }{2\left( \wp(z)+\wp (w) \right)\left( \wp(z)-\wp(w) \right)^2}-\frac{\wp(z)\wp(w)}{\wp(z)+\wp(w)}-\frac{g_2}{4\left( \wp(z)+\wp(w) \right)}\\
    \nonumber &=\frac{1}{4\wp(z)\wp(w)}\left( \frac{\wp(z)\wp'(w)-\wp(w)\wp'(z)}{\wp(z)-\wp(w)} \right)^2+\frac{g_3}{4\wp(z)\wp(w)}.
    \end{align}
\end{corollary}
Taking the limit $w\to z$, we obtain the classic and two new duplication formulas
\begin{align*}
    \wp(2z)&=-\frac{\wp''(z)\left( \wp'(z)^2-\wp(z)\wp''(z) \right)}{4\wp(z)\wp'(z)^2}-\frac{\wp(z)}{2}-\frac{g_2}{8\wp(z)}\\
    &=\frac{1}{4}\left( \frac{\wp'(z)^2-\wp(z)\wp''(z)}{\wp(z)\wp'(z)} \right)^2+\frac{g_3}{4\wp(z)^2},
\end{align*}
Above formula can be simplified further and written only in terms of $\wp(z)$ using $\left(\wp'\right)^2=4\wp^3-g_2\wp-g_3$ and $\wp''=6\wp^2-g_2/2$. Indeed, we have
\begin{align*}
    \wp(2z)&=\frac{16\wp(z)^4+8g_2\wp(z)^2+32g_3\wp(z)+g_2^2}{16\left( 4\wp(z)^3-g_2\wp(z)-g_3 \right)}\\
    &=\frac{1}{16} \frac{\left(4\wp(z)^3+g_2\wp(z)+2g_3\right)^2}{\wp(z)^2\left(4\wp(z)^3-g_2\wp(z)-g_3\right)}+\frac{g_3}{4\wp(z)^2}.
\end{align*}
Comparing the equations in Corollary \ref{ClassicalCorollary}, we have the following identities for the Weierstrassian elliptic invariants $g_2,g_3$. 
\begin{corollary}\label{InvariantCor}
    The following identity holds,
     \begin{align*}
        g_2&=4\left( \wp(z)+\wp(w) \right)^2-4\wp(z)\wp(w)-\left( \wp(z)+\wp(w) \right)\left( \frac{\wp'(z)-\wp'(w)}{\wp(z)-\wp(w)} \right)^2\\
        &-2\frac{\left( \wp'(z)-\wp'(w) \right)\left( \wp(z)\wp'(w)-\wp(w)\wp'(z) \right) }{\left( \wp(z)-\wp(w) \right)^2}.
     \end{align*}
     Furthermore, we have
     \begin{align*}
         g_3=\wp(z)\wp(w)\left( \frac{\wp'(z)-\wp'(w)}{\wp(z)-\wp(w)} \right)^2-4\wp(z)\wp(w)\left( \wp(z)+\wp(w) \right)-\left( \frac{\wp(z)\wp'(w)-\wp(w)\wp'(z)}{\wp(z)-\wp(w)} \right)^2.
     \end{align*}
\end{corollary}
Upon taking $z=\omega_1$ and $w=\omega_2$ in Corollary \ref{InvariantCor} we get the classic formulae $g_2=-4(e_1e_2+e_2e_3+e_3e_1)$ and $g_3=4e_1e_2e_3$. 
\end{example}

\begin{example}
    Let us take $n_1=2$, $\gamma_1=1$, $(k_1,k_2,k_3)=(1,0,-2)$ and $\ell=3$. We now note several novel three-term addition formulas. We have the system
\begin{gather*}
    \wp''(u)=\lambda_1\wp'(u)+\lambda_2\wp(u)+\lambda_3\\
    \wp''(v)=\lambda_1\wp'(v)+\lambda_2\wp(v)+\lambda_3\\
    \wp''(w)=\lambda_1\wp'(w)+\lambda_2\wp(w)+\lambda_3.
\end{gather*}
Using Frobenius-Stickelberger formula \eqref{FrobeniusStick}, condition \eqref{DetCondition} is equivalent to
\begin{align*}
    0\ne\begin{vmatrix}
            1 & \wp(u) & \wp'(u)  \\
            1 & \wp(v) & \wp'(v)  \\
            1 & \wp(w) & \wp'(w)
    \end{vmatrix}=2\frac{\sigma(u+v+w)\sigma(u-v)\sigma(v-w)\sigma(w-u)}{\sigma(u)^3\sigma(v)^3\sigma(w)^3}.
\end{align*}
This directly translates to the assumptions that $u,v,w$ are distinct $\bmod \Lambda_\Omega$ and $u+v+w\notin\Lambda_\Omega$. These conditions ensure that we do not encounter the poles of the $\wp$ function. The determinant in Theorem \ref{DetThm} unfortunately vanishes identically by the Frobenius-Stickelberger formula. We get that $\varphi(s)=\lambda_1^2\wp'(s)^2-\left( \wp''(s)-\lambda_2\wp(s)-\lambda_3 \right)^2$. Utilizing \eqref{WDE} and its corollary $\wp''=6\wp^2-g_2/2$, we calculate $\mu(r)$ as follows,
\begin{align*}
    &\mu(0)=-\lambda_1^2g_3-\lambda_3^2-\lambda_3g_2-\frac{g_2^2}{4}\\
    &\mu(1)=-\lambda_1^2g_2-2\lambda_2\lambda_3-\lambda_2g_2\\
    &\mu(2)=-\lambda_2^2+12\lambda_3+6g_2\\
    &\mu(3)= 4\lambda_1^2+12\lambda_2\\
    &\mu(4)=-36.
\end{align*}
By Cramer's rule, $\lambda_1,\lambda_2$ and $\lambda_3$ are given by
\begin{align*}
    &\lambda_1=\frac{\sum_{cyc} \wp''(u)\wp(v)-\wp(u)\wp''(v)}{\sum_{cyc} \wp'(u)\wp(v)-\wp(u)\wp'(v)}\\
    &\lambda_2=\frac{\sum_{cyc} \wp'(u)\wp''(v)-\wp''(u)\wp'(v)}{\sum_{cyc} \wp'(u)\wp(v)-\wp(u)\wp'(v)}\\
    &\lambda_3=\frac{\sum_{cyc} \wp(u)(\wp''(v)\wp'(w)-\wp'(v)\wp''(w))}{\sum_{cyc} \wp'(u)\wp(v)-\wp(u)\wp'(v)},
\end{align*}
where summations are cyclic in terms of $u,v,w$. Thus we have the following addition theorems.

\begin{corollary}
    We have the following addition formulae,
    \begin{align}\label{3Termaddition}
        \wp(u+v+w)&=\frac{\lambda_1^2}{9}+\frac{\lambda_2}{3}-\wp(u)-\wp(v)-\wp(w)\\
    \nonumber &=\frac{\lambda_2^2-12\lambda_3-6g_2}{36(\wp(u)+\wp(v)+\wp(w))}-\frac{\wp(u)\wp(v)+\wp(v)\wp(w)+\wp(w)\wp(u)}{\wp(u)+\wp(v)+\wp(w)}\\
    \nonumber &=-\frac{\lambda_1^2g_2+2\lambda_2\lambda_3+\lambda_2g_2}{36(\wp(u)\wp(v)+\wp(v)\wp(w)+\wp(w)\wp(u))}-\frac{\wp(u)\wp(v)\wp(w)}{\wp(u)\wp(v)+\wp(v)\wp(w)+\wp(w)\wp(u)}\\
    \nonumber &=\frac{4\lambda_1^2g_3+4\lambda_3^2+4\lambda_3g_2+g_2^2}{144\wp(u)\wp(v)\wp(w)}.
    \end{align}
\end{corollary}
Similarly, comparing above equations for different sides of $\wp(u+v+w)$ and rearranging gives us a variety of novel identities involving the invariants $g_2$ and $g_3$. Now we take the limits $v,w\to u=z$ to obtain triplication formula for $\wp(z)$. After a lengthy computation of limits we obtain
\begin{align*}
    &\lim_{u,v,w \to z} \lambda_1=\frac{\wp''(z)\wp'''(z)-\wp'(z)\wp^{(4)}(z)}{\wp''(z)^2-\wp'(z)\wp'''(z)}\\
    &\lim_{u,v,w \to z} \lambda_2=\frac{\wp''(z)\wp^{(4)}(z)-\wp'''(z)^2}{\wp''(z)^2-\wp'(z)\wp'''(z)}\\
    &\lim_{u,v,w \to z} \lambda_3=\frac{\wp(z)\wp'''(z)^2+ \wp'(z)^2\wp^{(4)}(z)+\wp''(z)^3-\wp(z)\wp''(z)\wp^{(4)}(z)-2 \wp'(z)\wp''(z)\wp'''(z)}{\wp''(z)^2-\wp'(z)\wp'''(z)}\\
\end{align*}
Hence we have the following triplication formula,
\begin{align}\label{triplication}
    \wp(3z)=\frac{1}{9}\left( \frac{\wp''(z)\wp'''(z)-\wp'(z)\wp^{(4)}(z)}{\wp''(z)^2-\wp'(z)\wp'''(z)} \right)^2+\frac{1}{3}\frac{\wp''(z)\wp^{(4)}(z)-\wp'''(z)^2}{\wp''(z)^2-\wp'(z)\wp'''(z)}-3\wp(z).
\end{align}
We furthermore get the more complicated ones
\begin{align*}
    \wp(3z)&=\frac{1}{108\wp(z)}\left( \frac{\wp''(z)\wp'''(z)-\wp'(z)\wp^{(4)}(z)}{\wp''(z)^2-\wp'(z)\wp'''(z)} \right)^2-\frac{g_2}{18\wp(z)}-\wp(z)\\
    &-\frac{1}{9\wp(z)}\frac{\wp(z)\wp'''(z)^2+ \wp'(z)^2\wp^{(4)}(z)+\wp''(z)^3-\wp(z)\wp''(z)\wp^{(4)}(z)-2 \wp'(z)\wp''(z)\wp'''(z)}{\wp''(z)^2-\wp'(z)\wp'''(z)}
\end{align*}
and
\begin{align*}
    \wp(3z)&=-\frac{\wp (z)}{3}-\frac{g_2}{108\wp(z)^2}\left( \frac{\wp''(z)\wp'''(z)-\wp'(z)\wp^{(4)}(z)}{\wp''(z)^2-\wp'(z)\wp'''(z)} \right)^2-\frac{1}{108\wp(z)^2}\left( \frac{\wp''(z)\wp^{(4)}(z)-\wp'''(z)^2}{\wp''(z)^2-\wp'(z)\wp'''(z)} \right)\\
    &\times\left(g_2+2\frac{\wp(z)\wp'''(z)^2+ \wp'(z)^2\wp^{(4)}(z)+\wp''(z)^3-\wp(z)\wp''(z)\wp^{(4)}(z)-2 \wp'(z)\wp''(z)\wp'''(z)}{\wp''(z)^2-\wp'(z)\wp'''(z)}\right)
\end{align*}
and
\begin{align*}
    \wp(3z)&=\frac{g_2^2}{144\wp(z)^3}+\frac{g_3}{36\wp(z)^3}\left( \frac{\wp''(z)\wp'''(z)-\wp'(z)\wp^{(4)}(z)}{\wp''(z)^2-\wp'(z)\wp'''(z)} \right)^2\\
    &+\frac{g_2}{36\wp(z)^3}\frac{\wp(z)\wp'''(z)^2+ \wp'(z)^2\wp^{(4)}(z)+\wp''(z)^3-\wp(z)\wp''(z)\wp^{(4)}(z)-2 \wp'(z)\wp''(z)\wp'''(z)}{\wp''(z)^2-\wp'(z)\wp'''(z)}\\
    &+\frac{1}{36\wp(z)^3}\left( \frac{\wp(z)\wp'''(z)^2+ \wp'(z)^2\wp^{(4)}(z)+\wp''(z)^3-\wp(z)\wp''(z)\wp^{(4)}(z)-2 \wp'(z)\wp''(z)\wp'''(z)}{\wp''(z)^2-\wp'(z)\wp'''(z)} \right)^2.
\end{align*}
We obviously recommend the formula \eqref{triplication} for practical use. Above equations can be written only in terms of $\wp(z)$ utilizing the following
\begin{align*}
    \wp'''=12\wp'\wp\qquad \text{and}\qquad \wp^{(4)}=120\wp^3-18g_2\wp-12g_3.
\end{align*}
Substituting these in the definitions of $\lambda_1,\lambda_2$ and $\lambda_3$ we obtain 
\begin{align*}
    &\lim_{u,v,w \to z} \lambda_1=\frac{48\wp'\left( 4\wp^3-g_2\wp-g_3 \right)}{48\wp^4-24g_2\wp^2-48g_3\wp-g_2^2}\\
    &\lim_{u,v,w \to z} \lambda_2=-\frac{ 576\wp^5 -96g_2\wp^3+288g_3\wp^2+36g_2^2\wp+24g_2g_3}{48\wp^4-24g_2\wp^2-48g_3\wp-g_2^2}\\
    &\lim_{u,v,w \to z} \lambda_3=\frac{192\wp^6+240g_2\wp^4+768g_3\wp^3-12g_2^2\wp^2-96g_2g_3\wp-96g_3^2+g_2^3}{96\wp^4-48g_2\wp^2-96g_3\wp-2g_2^2},
\end{align*}
where $\wp=\wp(z)$. For example equation \eqref{triplication} is equivalent to
\begin{align*}
    \left( 48\wp^4-24g_2\wp^2-48g_3\wp-g_2^2 \right)^2\wp(3z)&= 256\wp^9+768g_2\wp^7+6144g_3\wp^6+480g_2^2\wp^5-384g_2g_3\wp^4
    \\
    &+\left(768g_3^2-144g_2^3 \right)\wp^3-192g_2^2g_3\wp^2+\left(9g_2^4-384g_2g_3^2  \right)\wp\\
    &+8g_2^3g_3-256g_3^3
\end{align*}
\end{example}

We remark that one might obtain simpler formulas such as \eqref{ClassicAdditionEq} and \eqref{3Termaddition} for $\wp(z_1+\cdots+z_\ell)$ by considering the system
\begin{align*}
    \wp^{(\ell-1)}(z_j)=\lambda_1\wp^{(\ell-2)}(z_j)+\lambda_2\wp^{(\ell-3)}(z_j)+\cdots+\lambda_{\ell-1}\wp(z_j)+\lambda_\ell.
\end{align*}
This corresponds to the case $n_1=\ell-1$, $\gamma_1=1$, $(k_1,\cdots, k_{\ell-1},k_\ell)=(\ell-2,\cdots ,0,-2)$. Similarly, the determinant in condition \eqref{DetCondition} can be calculated by Frobenius-Stickelberger formula and is equivalent to the statement that $z_i$ are distinct $\bmod\Lambda_\Omega$ with $z\notin\Lambda$. Taking the limits $z_1,\cdots,z_\ell\to s$, we get the $\ell$-multiplication formula for $\wp(s)$. \\

\noindent \textbf{Disclosure Statement.} The author(s) report no potential conflict of interest.

\end{document}